\documentclass[12pt,oneside,a4paper,leqno]{article}
\usepackage{psfrag}
\usepackage{subfigure}

\usepackage[all]{xy}
\CompileMatrices
\usepackage{geometry}
\usepackage{graphicx}
\usepackage{amsfonts,amssymb}
\usepackage[centertags]{amsmath}

\usepackage{amsthm}
\newcounter{theorems}

\theoremstyle{plain}
\newtheorem{theoAA}[theorems]{Theorem}

\swapnumbers
\newcounter{lemma}
\numberwithin{equation}{section}

\newtheoremstyle{par}%
     {\topsep}%
     {\topsep}%
     {\itshape}%
     {}%
     {\bfseries}%
     {}%
     {.5em}%
     {}%

\newtheoremstyle{parrm}%
     {\topsep}%
     {\topsep}%
     {\normalfont}%
     {}%
     {\itshape}%
     {}%
     {.5em}%
     {}%
\theoremstyle{plain}
\numberwithin{equation}{section}

\newtheorem{propo}[equation]{Proposition}

\theoremstyle{definition}
\newtheorem{defi}[equation]{Definition}

\theoremstyle{remark}
\newtheorem{remark}[equation]{Remark}

\theoremstyle{par}

\newtheorem{lemma}[equation]{}

\theoremstyle{parrm}
\newtheorem{rem}[equation]{}

\makeatletter
\def\tagform@#1{\maketag@@@{\ignorespaces#1\unskip\@@italiccorr}}

\makeatother
\newcommand{\CC}{\mathbb{C}}
\newcommand{\R}{\mathbb{R}}
\newcommand{\ze}{\mathbb{Z}}
\newcommand{\ee}{\mathbf{e}}
\newcommand{\st}{\ \mathrm{|} \ }
\newcommand{\from}{\colon}
\newcommand{\X}{\mathcal{X}}
\newcommand{\bfomega}{\mathbf{\underline\omega}}
\newcommand{\T}{\mathbb{T}}
\newcommand{\action}{\mathcal{A}}

\newcommand{\nn}{\mathbf{n}}
\newcommand{\ff}{\varphi}
\newcommand{\dt}{\, dt}

\theoremstyle{par}

\theoremstyle{plain}

\theoremstyle{definition}

\begin{document}
\pagenumbering{arabic}
\title{%
Symmetry groups and non-planar collisionless action-minimizing solutions
of the three-body problem in three-dimensional space}

\author{
Davide L.~Ferrario\footnote{%
Partially supported by MIUR, National Project 
``Variational Methods and Nonlinear Differential Equations''.
Address: 
Dipartimento di Matematica e Applicazioni, University of Milano--Bicocca,
Via Cozzi 53, 20125 Milano (Italy).
e-mail: \emph{ferrario@matapp.unimib.it}.
}
}
\date{\today}
\maketitle
\begin{abstract} 

Periodic and quasi-periodic solutions of the $n$-body
problem  can be found as minimizers of the
Lagrangian action functional restricted to suitable spaces of symmetric 
paths. The main purpose of this paper is to develop
a systematic approach to the equivariant minimization 
for the three-body problem in the three-dimensional space.
First 
we give a finite complete list of symmetry groups fit to the minimization of
the action, with the property that any other 
symmetry group can be reduced to be
isomorphic to one of these representatives.  
A second step is to  prove that
the resulting (local and global) symmetric action-minimizers are always
collisionless (when they are not already bound to collisions).
Furthermore, we  prove some results  addressed 
to the question whether minimizers are planar 
or non-planar; as a consequence of the theory  
we will give  general criteria for a 
symmetry group to  yield planar or homographic minimizers  
(either homographic or not, as in the Chenciner-Montgomery eight
solution);
on the other hand 
we will provide a rigorous proof of 
the existence of some interesting one-parameter
families of periodic and quasi-periodic non-planar orbits.  
These include the choreographic
Marchal's $P_{12}$ family with equal masses -- together with a less-symmetric
choreographic family (which anyway probably coincides with the $P_{12}$).

\vspace{12pt}
\noindent
\emph{MSC Subj. Class}: 
Primary 
70F10;
Secondary 
70G75,
37C80,
70F16,
37J45.

\vspace{12pt}
\end{abstract}
\section{Introduction}
\label{sec:intro}
In some recent  papers classical variational
methods have been successfully
applied in the proof of the existence of periodic or
quasi-periodic solutions for  the $n$-body problem.
Suitable symmetry groups of the Lagrangian action functional
have been introduced and exploited in order to apply
the aforementioned techniques to the class
of all symmetric loops, to name a few, in the 
articles 
\cite{%
bessi,
chen,%
MR2032484,%
monchen,%
chenven,%
FT2003,%
marchal,%
montgomery_contmath,%
serter1}.
Surveys and further details 
on this approach
can be found for example in 
\cite{%
zz,
chenICM,chenciner_taiyuan,
FT2003,
montgomery_notices%
}.
The major problem in the search of equivariant minimizers is that 
of collisions: a (local or global) minimizer might consist
\emph{a priori} of a colliding trajectory. 
The latest significant breakthrough in this direction 
has been allowed by  Marchal's averaging technique 
\cite{chenICM,FT2003}. In the paper \cite{FT2003} the authors develop
a general theory for $G$-equivariant minimizers
and present a class of groups that yield always,
as a consequence of Marchal's 
averaging technique, collision-free minimizers.
In \cite{zz} this result was extended to all possible symmetry
groups for the planar three-body problem.
A naturally related problem is 
to find and classify all possible symmetry groups
and to understand whether the resulting minima
are  rotating central configurations 
or if they are new solutions (and, at the same time,
to provide rigorous proofs of the existence and of the properties
of some solutions whose existence was 
accepted as a fact after numerical evidence). 
This has been done for the planar problem 
 by V. Barutello, S. Terracini and the author
in \cite{zz}.
The purpose of the paper is to give a complete
answer to the classification problem for the three-body
problem in the space 
and at the same time to determine and describe
properties of the resulting minimizers. 
In particular,
we focus on non-planar orbits,
since planar orbits have been already included in the list of
 \cite{zz}.
In order to state the main results, we anticipatively sketch some basic
definitions:
A symmetry group $G$ of the Lagrangian
functional $\action$ (see \ref{def:sg} below)
is termed \emph{bound to collisions}
if all $G$-equivariant loops actually have collisions
(see \ref{def:btc} below),
\emph{fully uncoercive} 
if for every possible rotation vector $\bfomega$ 
the action functional 
$\action_\bfomega^G$ 
in the frame rotating around $\bfomega$ with angular
speed $|\bfomega|$
is not coercive in the space of $G$-equivariant loops (that is,
its global minimum escapes to infinity -- see \ref{defi:fullyuncoercive});
moreover, $G$ is termed 
\emph{homographic} if all $G$-equivariant loops
are constant up to orthogonal motions and rescaling.
Note that if there is a rotation axis $\bfomega$ then
the group $G$ is implicitly assumed to be a symmetry 
group of the action functional $\action_\bfomega$ in the rotating
frame (that is, the functional including the centrifuge
and Coriolis terms); such a group is termed \emph{of type R}
(see \ref{defi:typeR3d} below); 
finally, 
the \emph{core} of the group $G$ is the subgroup
of all the elements which do not move the time $t\in \T$
(see \ref{defi:core} below).
In the first theorem we classify symmetry groups,
up to change in rotating frame. For the symbols used
we refer to sections \ref{sec:prelim} and 
\ref{sec:onedim} below.
\begin{theoAA}
\label{theo:mt1}
Symmetry groups not bound to collisions, not
fully uncoercive and not homographic are, up to a change
of rotating frame,   either
the three-dimensional extensions of planar groups (if trivial core)
listed in table 
\ref{table:thisclass} or 
the vertical isosceles triangle \ref{defi:vip}  (if non-trivial core).
\end{theoAA}

\begin{table}[h!]
\caption{Space extensions of planar 
symmetry groups with trivial core}
\label{table:thisclass}
\begin{center}
\begin{tabular}{l|r}
\hline
\emph{Name} & \emph{Extensions}  \\
\hline
Trivial             & $C_1^-$\\
Line                & 
$L_2^{+,-}$, 
$L_2^{-,+}$
 \\
Isosceles           & 
$H_2^{+,-}$,
$H_2^{-,+}$
\\
Hill                & 
$H_4^{+,-}$,
$H_4^{-,+}$
\\
$3$-choreography    & 
$C_3^+$, $C_3^{-}$ \\
Lagrange            & 
$L_6^{+,+}$,
$L_6^{+,-}$,  
$L_6^{-,+}$
\\
$D_6$               & 
$D_6^{+,-}$,
$D_6^{-,+}$
 \\
$D_{12}$            & 
$D_{12}^{-,+}$ 
\\
\end{tabular}
\end{center}
\end{table}
The next theorem is the answer to the natural questions about
collisions and description of some main features of minimizers.
\begin{theoAA}
\label{theo:mt2}
Let $G$ be a symmetry group not bound to collisions and 
not fully uncoercive. Then
\begin{enumerate}
\item Local minima of $\action_\bfomega^G$ do not have collisions.
\item In the following cases minimizers are planar trajectories:
\begin{enumerate}
\item If  $G$ is not of type R: $D_6^{+,-}$,
$D_6^{-,+}$ and $D_{12}^{-,+}$
(and then $G$-equivariant minimizers are 
Chenciner--Montgomery eights).
\item If there is a $G$-equivariant minimal Lagrange rotating solution:
$C_1^-$, 
$H_2^{+,-}$,
$C_3^+$,
$L_6^{+,+}$,
$L_6^{+,-}$ 
(and then the Lagrange solution is of course the minimizer).
\item
If the core is non-trivial and 
it is not the vertical isosceles \ref{defi:vip}
(and then minimizers are homographic).
\end{enumerate}
\item  
In the following cases minimizers are always non-planar:
\begin{enumerate}
\item 
The groups $L_6^{-,+}$
and $C_{3}^-$ for all $\omega\in (-1,1) + 6\ze$, $\omega\neq 0$
(the minimizers for $L_6^{-,+}$  are the elements of 
 Marchal family $P_{12}$, and minimizers 
of $C_3^-$ are  a less-symmetric family $P_{12}'$
\footnote{Highly likely they are not distinct families: this is 
the recurring phenomenon of ``more symmetries than expected''
in $n$-body problems.}).
\item 
The extensions of line and Hill-Euler type groups,
for on open subset 
  of mass distributions and angular speeds $\omega$
(explicitely given in \ref{eq:othervar2}):
$L_2^{+,-}$,
$L_2^{-,+}$,
$H_4^{+,-}$
and $H_4^{-,+}$ (for $L_2^{-,+}$
this happens also with equal masses).
\item The vertical isosceles \ref{defi:vip} for 
suitable choices of masses and $\omega$.
\end{enumerate}
\end{enumerate}
\end{theoAA}
In the article we develop the needed tools
and prove  these statements, after the necessary
explanations about preliminary results and notation. Together
with the results of \cite{FT2003} and \cite{zz} it  is 
exposed a theory of action-minimizing 
symmetry periodic $n$-body orbits. 
In section \ref{sec:prelim} we introduce all the 
definitions needed in the sequel and prove some
preliminary results. In section \ref{sec:onedim} we introduce
the concept of three-dimensional \emph{extension} 
of a planar symmetry group so that we can use 
the classification of \cite{zz}. In section \ref{sec:nottyper}
the angular momentum $J$ enters into account,
and we show how the existence of rotation axes
is related to the possibility of being non-zero of $J$.
Afterward, in section \ref{sec:vertical}
we prove some interesting estimates on second variations,
which are remarkably simple (incidentally, they
work not only for $3$ bodies, but for $n$ arbitrary).
It is by an application of these
simple estimates that one proves the fact that the non-planar
quasi-periodic orbits listed in theorem \ref{theo:mt2} exist.
In section \ref{sec:zzz} we come to the classification
of three-dimensional space symmetries, which is 
a proof of theorem \ref{theo:mt1}.
The proof of the various items of theorem \ref{theo:mt2}
is done in  section \ref{sec:main}.
Finally, in section \ref{sec:remarks} some concluding remarks
are collected.

Before we start with the next section, a few words have to be spent
on the existence of the $P_{12}$.
A different -- and very elementar --
proof of the existence of the $P_{12}$-family
with $D_{12}$-symmetries was presented by A. Chenciner 
in \cite{chenICM,chenciner_taiyuan}, which does not 
require local results on collisions, since collisions
are excluded by action level estimates. The advantage 
of our approach is that it can be plainly extended to the case
of any odd number $n\geq 3$ of bodies in the space 
(see remark \ref{rem:extendit} below).
All other results are, to our knowledge, new: 
whenever similar methods or results were published elsewhere, 
it has been remarked in-place.
\section{Preliminaries}
\label{sec:prelim}

Consider the linear space of configurations
with center of mass in $0$
\[
\X  = \{ x=(x_1,x_2,x_3) \in E^3 \st m_1x_1 + m_2x_2 + m_3x_3=0 \}.
\]
Let $\T=S^1 \cong \R/2\pi\ze$
be  the unit circle of length $2\pi$. 
We are dealing with periodic orbits of the Newtonian $n$-body
problem, which will be seen as critical points of a suitable
functional on the Sobolev space 
$\Lambda = H^1(\T,\X)$ 
 consisting of all 
$L^2$ loops $\T\rightarrow\X$ with $L^2$
derivative. It is an Hilbert space with the scalar product
\[
x \cdot y = \int_{\T}(x(t)y(t)+\dot{x}(t)\dot{y}(t))dt.
\]

The $\alpha$-homogeneous
Newtonian potential can be written as
\begin{equation}
\label{pot:function}
U(x)= \frac{m_1 m_2}{|x_1-x_2|^\alpha}+
\frac{m_1 m_3}{|x_1-x_3|^\alpha}+\frac{m_2 m_3}{|x_2-x_3|^\alpha}.
\end{equation}
Let $\bfomega\in E\cong \R^3$ be a vector.  The kinetic form 
in a frame uniformly rotating around $\bfomega$
with angular speed $\omega=|\bfomega|$ is  defined by
\begin{eqnarray}
\label{kin:energy}
2 K(x,\dot{x})&=&
\sum_{i=1}^3
m_i | \dot x_i + \Omega x_i |^2,
\end{eqnarray}
where $\Omega$ is the matrix 
\[
\left(
\begin{array}{ccc}
0 & \omega_3 & -\omega_2 \\
-\omega_3 & 0 & \omega_1 \\
\omega_2 & -\omega_1 & 0 \\
\end{array}
\right)
\]
obtained by the coefficients $(\omega_1,\omega_2,\omega_3)$
of $\bfomega \in \R^3$.
Thus, for every $\bfomega$ the Lagrangian can be written as
\begin{equation}
\label{lagrangian}
L_\bfomega(x,\dot{x})=L_\bfomega=K_\bfomega+U,
\end{equation}
and, finally, the action functional as
\begin{equation}
\label{action:functional}
\action_\bfomega(x)=
\int_{\T}L_\bfomega(x(t),\dot{x}(t)) dt.
\end{equation}

\begin{defi}\label{def:sg}
We term 
\emph{symmetry group} every subgroup of 
$O(\T) \times O(3) \times \Sigma_3$,
where $O(\T) = O(2)$ is the orthogonal group of dimension $2$
acting on the time circle,
$O(3)$ the orthogonal group of dimension $3$ acting on the 
space $E$ and $\Sigma_3$ the symmetric group on the three
elements $\{1,2,3\}$. 
\end{defi}

Given a subgroup $G\subset O(\T) \times O(3) \times \Sigma_3$,
it is possible to define three homomorphisms
$\tau\from G \to O(\T)$,
$\rho\from G \to O(3)$ 
and $\sigma\from G \to \Sigma_3$ 
by projection onto the first, second or third
factor of the direct product.
Given $\tau$, $\rho$ and $\sigma$
one can  define in an obvious way
a $G$-action on $\T$, $E$ and $\{1,2,3\}$,
and hence an action  on the centered configuration space  $\X$,
provided that for every $g\in G$ it happens that 
$m_{i} = m_{gi}$ (that is, masses are constant in 
$G$-orbits in the index set $\{1,2,3\}$.
Thus, there is an  induced action of $G$ 
on the Sobolev space $\Lambda$ of loops, 
defined by $g (x(t)) = (gx)(g^{-1}t)$. The action is orthogonal
on $\Lambda$ so that Palais theorem says that, for a given 
group $G$, 
if the functional $\action_\bfomega\from \Lambda \to \R$ 
is $G$-invariant 
and a $G$-equivariant  loop $x(t)$  is 
 collisionless and  critical  for the restriction  $\action_\bfomega^G = 
\action_\bfomega |\Lambda^G$,
then  $x(t)$ is critical for $\action_\bfomega$.
\begin{defi}\label{def:btc}
A group is termed \emph{bound to collisions}
if for  every equivariant loop $x(t) \in \Lambda^G$
collisions occur, that is, 
for each $G$-equivariant $x(t)$ 
there is $t_c\in \T$ and $i\neq j \in \nn=\{1,2,3\}$
such  that $x_i(t_c) = x_j(t_c)$.
\end{defi}

\begin{defi}\label{defi:homographic}
A group $G$ is termed \emph{homographic}
if every equivariant loop $x(t) \in \Lambda^G$ 
is constant up to rescaling and orthogonal motions.
\end{defi}

\begin{defi}\label{defi:core}
The kernel $\ker\tau$ is termed the \emph{core} of the symmetry
group $G$.
\end{defi}

\begin{defi}\label{cyclic}
A group $G$ is termed of \emph{cyclic}, \emph{brake}
or \emph{dihedral} type respectively
if $G/\ker\tau$  acts orientation-preserving 
on the time circle $\T$, 
if 
$G/\ker\tau$ has order $2$ and acts
orientation-reversing on $\T$
or if $G/\ker\tau$ is a dihedral group of order $\geq 4$.
\end{defi}
Consider the following elements in $O(2)$:
$1$ is the trivial motion,
$-1$ is the rotation of angle $\pi$ and  $l$ 
is a    reflection along a line.
Elements of the symmetric group $\Sigma_3$ will
be denoted in the cyclic permutation notion.
Let $\ker\det(\tau)\subset G$ denote the 
subgroup of $G$ of the elements acting orientation-preserving
on $\T$.
A symmetry group with trivial core will be fully determined 
once the images $\rho(r)$ and $\sigma(r)$ 
of a generator $r$ 
of  $\ker\det(\tau)$ 
rotating a minimal angle are given,
together, if it is not of cyclic type, with 
the images $\rho(h)$ and $\sigma(h)$ 
of one of the elements not in $\ker\det(\tau)$
(which have order $2$).
The full list of representatives of the planar 
classification exposed in \cite{zz} can be therefore
found in table \ref{table:class}, with 
the corresponding generators.
\begin{table}[ht]
\caption{Planar symmetry groups with trivial core}
\label{table:class}
\begin{center}
\begin{tabular}{l|c|cc}
\hline
\emph{Name} & \emph{Symbol}  & $\rho(r)$,$\sigma(r)$ & 
$\rho(h)$,$\sigma(h)$\\
\hline
Trivial             & $C_1$ & 1,()  & \\
Line                & $L_2$ & 1,() & l,() \\
$2$-$1$-choreography & $C_2$ & 1, (1,2) &  \\
Isosceles           & $H_2$ & 1,() & l,(1,2) \\
Hill                & $H_4$ & 1,(1,2) & l,(1,2) \\
$3$-choreography    & $C_3$ & 1,(1,2,3) & \\
Lagrange            & $L_6$  & 1,(1,2,3) & l,(1,2) \\
$C_6$               &  $C_6$  & l,(1,2,3) & \\
$D_6$               & $D_6$ & 1,(1,2,3) & -1,(1,2) \\
$D_{12}$            & $D_{12}$ l,(1,2,3) & -1,(1,2) \\
\end{tabular}
\end{center}
\end{table}
\begin{defi}
\label{defi:typeR}
A planar symmetry group $G$ is said \emph{of type R} 
if the determinant
homomorphisms $\det(\rho), \det(\tau) \from G \to \{+1,-1\}$
coincide, that is, if they coincide as 
$G$-representations.
\end{defi}
\begin{defi}
\label{defi:typeRdir}
A vector $v\in E\cong \R^3$ 
is termed a \emph{rotation axis}  with respect
to a symmetry group $G$ 
if the line spanned by $v$ in $E$ is $G$-invariant  and 
the following equality of one-dimensional 
$G$-representations holds:
\[
\det (\tau)  \det (\rho) = \det (v),
\]
where $\det(v)$ denotes the  real representation
of $G$ induced by restricting $\rho$ to the invariant subspace
generated by $v\in E$.
\end{defi}
\begin{lemma}
\label{typeR3d}
The restriction of a three-dimensional symmetry group $G$  to 
the orthogonal complement of a rotation axis $v\in E$
is a planar symmetry group of type R.
Conversely, if the restriction of $G$ to an invariant plane  
is of type R, then the orthogonal complement of the invariant 
plane is a rotation axis for $G$.
\end{lemma}
\begin{proof}
Let $\tau$, $\rho$ and $\sigma$ be the defining
homomorphisms of $G$, where $\rho\from G \to O(3)$ 
can be written as $\rho=\rho_2 \times \rho_1$, 
with 
$\rho_2\from G \to O(2)$ 
induced by restriction to the (invariant)  orthogonal complement
of $v$ and 
$\rho_1\from G\to O(1)$ by restriction to $v$.
Since 
\begin{equation}\label{eq:brr}
\det (\rho) = \det (\rho_2) \det(\rho_1) =
\det(\rho_2) \det(v),
\end{equation}
the planar symmetry group defined by  $\tau,\rho_2,\sigma$ 
by \ref{defi:typeR} is of type R 
if and only if  
$\det(\rho_2) = \det(\tau)$,
and hence if an only if 
\( \det(\rho) \det(v) = \det(\tau) \) as claimed.
\end{proof}

The previous lemma yields the following natural definition.
\begin{defi}
\label{defi:typeR3d}
A space symmetry group $G$ is said 
\emph{of type R}
if it has at least one rotation axis (that is,
if it is the extension of a planar group of type R.
\end{defi}
\begin{lemma}
Let $\bfomega\in E\cong \R^3$ be a rotation axis for a symmetry
group $G$. Then the Lagrangian action functional  $\action_\bfomega$
(defined 
in \ref{action:functional})
in a frame  rotating around $\bfomega$ with angular speed $\omega=|\bfomega|$
is $G$-invariant.
\end{lemma}
\begin{proof}
Let $g\in G$ and   $x(t) \in \Lambda^G$.
Since for every $i\in \{1,2,3\}$
\[
x_{gi}(\tau(g)t) = \rho(g)x_i(t),
\]
the derivative fulfills the equality
\[
\det(\tau(g)) \dot x_{gi}(\tau(g)t) = \rho(g)\dot x_i(t).
\]
Thus for every $g\in G$ and $t\in \T$,
\[
\dot x_{gi}(\tau(g)t)  + \Omega x_{gi}(\tau(g)(t)) =
\det(\tau(g))\rho(g) \dot x_{i} + \Omega \rho(g) x_{i}(\tau(t)),
\]
and hence
\[
\left| \dot x_{gi}(\tau(g)t)  + \Omega x_{gi}(\tau(g)(t)) \right|^2 =
\left|\dot x_{i} + 
\det(\tau(g))\rho(g^{-1})\Omega \rho(g) x_{i}(\tau(t)) \right|^2.
\]
One can deduce that the action functional $\action_\bfomega$ 
is $G$-invariant if (and only if)
for every $g$ 
\begin{equation}\label{eq:this}
\det(\tau(g))\rho(g^{-1})\Omega \rho(g) = \Omega.
\end{equation}
If $\Omega\neq 0$,
equation \ref{eq:this} holds 
if and only if $\det(\rho_2) = \det(\tau)$, 
where as above $\rho_2$ denotes the restriction
of $\rho$ to the plane orthogonal to $\bfomega$. 
But by \ref{eq:brr} this is equivalent to the identity
$\det(\rho)\det(\bfomega) = \det(\tau)$,
that is,  $\bfomega$ is  a rotation axis
as in definition \ref{defi:typeRdir}.
\end{proof}
\begin{defi}
\label{defi:fullyuncoercive}
A symmetry group $G$ is said \emph{fully uncoercive} if
for every possible rotation vector $\bfomega$,
the action functional $\action^G_\bfomega$ is 
not coercive.
\end{defi}

The following proposition is an easy consequence of  the definition and 
(4.1) of \cite{FT2003}.
\begin{lemma}\label{lem:pre}
Let $G$ be a symmetry group.
\begin{enumerate}
\item 
If there are no rotation axes and $\X^G \neq 0$, 
(or, equivalently, $\action^G$ is not coercive),
then $G$ is fully uncoercive.
\item 
If every rotating axis is  uncoercive as a one-dimensional
$G$-module and the action on the index set is not transitive, 
then 
$G$ is fully uncoercive.
\end{enumerate}
\end{lemma}

\begin{defi}\label{defi:rotating}
If $G$ and $G'$ are two groups conjugated in
$O(\T)\times O(3) \times \Sigma_3$,
we will write $G\cong G'$.
If  there exists a change of rotating frame
for which 
a group $G$ can be written as $G'$, which is 
conjugate to a third group $G''$, 
we will write $G\sim G''$.
It is easy to see that $\cong$ and 
$\sim$ are equivalence relations, 
and that $G\cong G' \implies |G|=|G'|$,
while the same does not hold for $\sim$ (see \cite{zz}, section 3
for further details on changing the coordinates in a rotating frame).
\end{defi}
\begin{propo}\label{propo:gordon}
Let $G$ be a symmetry group such that a Lagrange
rotating solution $x(t) = \{x_j(t)\} = \{e^{ikt}\zeta_3^j\}$ is 
$G$-equivariant and  $|k+\omega|$ is minimal (as $k$ varies
in $\ze$) and not zero.
Then $x(t)$  is the absolute minimum of the action functional.
\end{propo}
\begin{proof}
This is, for the three-dimensional plane, proposition  (4.1) of
\cite{zz} (see also \cite{chendesol}).
Actually, for the three body problem the proof is
straightforward: assume that $\sum_{i} m_i = 1$;
since  the center of mass is in zero $\sum_{i} m_i x_i = 0$,
and hence $\sum_i m_i \dot x_i = 0$, the kinetic 
energy can be written in terms of the  differences
\begin{equation}\label{eq:differences}
\dfrac{1}{2} \sum_{i} m_i|\dot x_i + \Omega x_i|^2 = 
\dfrac{1}{2} 
\sum_{i<j} m_i m_j |\dot x_i  -\dot x_j + \Omega(x_i - x_j) |^2.
\end{equation}
Thence the action functional is written as the sum of 
three terms of the type
\[
\dfrac{1}{2}m_im_j|\dot x_i - \dot x_j + \Omega x_i - \Omega x_j|^2
+ m_i m_j |x_i - x_j|^{-\alpha},
\]
which is a Kepler (one-center) problem for the 
variable $y=x_i - x_j$ in the rotating frame.
Since a rotating solution $x(t)$ with $|k+\omega|$ minimal exists
by assumption, it 
yields three (identical, up to a time-shift) 
rotating solutions in $y$, with $|k+\omega|$ minimal.
It is easy to conclude the proof and to show 
that every trajectory has an action
which is at least three times the action of a minimal one-center $y$.
\end{proof}
\section{Three-dimensional extensions of planar symmetry groups}
\label{sec:onedim}
In this section we will take planar groups,
listed in table \ref{table:class},
and define some extensions acting on the three-dimensional
space. Of all the resulting groups, 
we will take into account only the extensions with trivial
core, not bound to collisions and not fully uncoercive. 
The outcome is the list of table \ref{table:thisclass}.
We proceed as follows.
Consider one of the planar groups in table \ref{table:class}.
It can be extended to a group acting on the three-dimensional space
simply by adding a one-dimensional real representation.
Now, since the groups have trivial core as 
for table \ref{table:class} we can assume that 
the symmetry group is generated by two elements
$r$ and $h$ with the following properties:
$\tau(r)$ is a time-shift in $\T$ of minimal angle
and $\tau(h)$ is a time-reflection (which exists only if the group
is not of cyclic type).
Up to conjugacy or change of orientation in $\T$ the
choice of $r$ and $h$ yields uniquely back the symmetry group
$G$.

Consider first the case of cyclic type, and let 
$r$ denote the cyclic generator  above
(in the notation of table \ref{table:class}, groups of cyclic type are
$C_1$, $C_2$, $C_3$ and $C_6$).
Now, $\rho_2(r) \in O(2)$ can be extended in two ways
to a matrix in $O(3)$: adding either a trivial
one-dimensional representation  or a non-trivial one.
Thus, for each cyclic group $C_i$ listed  above there exist two
 corresponding 
groups, denoted by $C_i^+$ and $C_i^-$,
which are generated by the element
$(\tau(r),\rho_2(r),\epsilon, \sigma(r))$
in 
$O(\T) \times O(2)\times O(1) \times \Sigma_3$,
for $\epsilon=\rho_1(r) \in\{\pm 1\}$.
The other cases can be dealt with in an analogous  way:
the choices 
are $2^2$: a sign for $\rho_1(r)$ and a sign for $\rho_1(h)$.
So, if $G$ is a symmetry group not of cyclic type,
its three-dimensional extensions groups will be denoted by
the symbol $G^{\epsilon_1,\epsilon_2}$, where 
$\epsilon_1$ is the sign of $\rho_1(r)$ and 
$\epsilon_2$ the sign of $\rho_1(h)$.
By \ref{typeR3d} 
the third axis will be a rotation axis
if and only if the planar 
symmetry group is of type R.
Furthermore, it is easy to see that 
if the action of 
the group on the index set is not transitive,
extensions of type
$+$, $+,+$ are fully uncoercive.
The list of remaining symmetry groups is therefore:
$C_1^-$, $C_1^{+,-}$, $C_1^{-,+}$, $C_1^{-,-}$, 
$L_2^{+,-}$, 
$L_2^{-,+}$, 
$L_2^{-,-}$
$C_2^-$  
$H_2^{+,-}$,
$H_2^{-,+}$,
$H_2^{-,-}$
$H_4^{+,-}$,
$H_4^{-,+}$,
$H_4^{-,-}$
$C_3^+$, $C_3^{-}$ 
$L_6^{+,+}$,
$L_6^{+,-}$,  
$L_6^{-,+}$,  
$L_6^{-,-}$,
$C_6^+$,
$C_6^-$
$D_6^{+,+}$,
$D_6^{+,-}$,
$D_6^{-,+}$,
$D_6^{-,-}$,
$D_{12}^{+,+}$ ,
$D_{12}^{+,-}$, 
$D_{12}^{-,+}$  and
$D_{12}^{-,-}$. 
Now, some of them are the the same after a change in coordinates:
$C_1^{-,-} \cong C_1^{-,+}$,
$L_2^{-,+} \cong L_2^{-,-}$,
$H_2^{-,+} \cong H_2^{-,-}$,
$H_4^{-,+} \cong H_4^{-,-}$,
$C_6^+ \cong C_3^-$,
$C_6^- \sim C_3$,
$D_6^{+,+} \cong L_6^{+,-}$,
$D_6^{-,+} \cong D_6^{-,-} \cong D_{12}^{+,-}$,
$L_6^{-,+} \cong L_6^{-,-} \cong D_{12}^{+,+}$ and 
$D_{12}^{-,-} \sim L_6^{+,-}$.
Furthermore, 
$C_1^{+,-}$, $C_1^{-,+}$ and 
$C_2^-$ are clearly fully uncoercive.
Hence the following lemma holds. 

\begin{lemma}
\label{lemma:onedim}
Of all the three-dimensional extensions of planar symmetry groups,
those  with trivial core, not bound to collisions 
and not fully uncoercive
are listed in table \ref{table:thisclass}.
\end{lemma}

\begin{remark}
The order of the space group now does not necessarily coincide with
the order of the planar group: for example, the order 
of $C_3^-$ is $6$ and not $3$. 
\end{remark}
The following lemma will be used as a key-step for 
the classification below.
\begin{lemma}
\label{lemma:sumofonedim}
Let $G$ a symmetry group for the three-body problem with trivial core. 
Then,
up to a change of rotating frame, 
$\rho$ is the sum of one-dimensional real representations.
\end{lemma}
\begin{proof}
Since $\ker\tau=1$, $G$ is  isomorphic to a subgroup 
of a finite dihedral group, and hence its orthogonal
irreducible representations  have dimension at most $2$.
So, $\rho$ can be written as $\rho_2\times \rho_1$ 
where $\rho_2\from G \to O(2)$ and $\rho_1\from G \to O(1)$.
Now, by (5.1) of \cite{zz} 
up to a change in rotating frame one can assume that 
$\rho_2(g)^2=1$ for every $g\in G$, so that 
$\rho_2$ is reducible as a sum of two one-dimensional
$G$-representations.   
\end{proof}
\section{Groups without rotation axes}
\label{sec:nottyper}

As we have seen in \ref{lemma:onedim},
the list of candidates for space symmetry groups
is given in table \ref{table:thisclass}.
Now we consider the $10$ groups yielded by extending
the three planar groups not of type R.

First, we prove a three-dimensional analogue of proposition
(3.9) of \cite{zz}.
Given a path $x(t) \in \Lambda$,
its angular momentum $J$ is the function of   $t\in \T$
given by
\begin{equation}
\label{eq:angmom}
J(t) =
\sum_{i\in \nn} m_i x_i \times \dot x_i 
\end{equation}
where $\times $ is the vector product in $E\cong \R^3$.
If $x$ is a (generalized) solution, then the angular momentum
is constant.
\begin{lemma}
\label{lemma:J}
For every equivariant $x(t) \in \Lambda^{G}$ and every $t\in \T$
\[
J(gt)  = \det(\rho(g)) \det(\tau(g)) \rho(g) J(t).
\]
\end{lemma}
\begin{proof}
It follows from the chain of equalities 
\begin{eqnarray*}
{J(\tau(g)t)} & =  &
\sum_{i=1}^3  m_i x_i(\tau(g)t) \times  \dot x_i(\tau(g)t)
\\
&=&
\sum_{i=1}^3 m_i \left[
(\rho(g) x_{g^{-1}i} (t)) \times
( \det(\tau(g)) \rho(g) \dot x_{g^{-1}i} (t)  )
\right] \\
&=& \det(\tau(g)) \sum_{i=1}^3 m_{g^{-1}i} \left[
(\rho(g) x_{g^{-1}i} (t)) \times
( \rho(g) \dot x_{g^{-1}i} (t)  )
\right] \\
&= &
\det(\tau(g)) \sum_{i=1}^3 m_{g^{-1}i}  \det(\rho(g)) \rho(g) \left[
x_{g^{-1}i} (t) \times
\dot x_{g^{-1}i} (t)
\right] \\
&=&
\det(\tau(g)) \det(\rho(g)) \rho(g) J(t). 
\end{eqnarray*}
\end{proof}
\begin{lemma}
\label{lemma:Estar}
Let $x\in \Lambda^G$ a $G$-equivariant periodic orbit with
angular momentum $J$.
Then 
$J$ belongs to the subspace   $E^*\subset E$
fixed by the $G$-representation 
$\det(\rho) \det(\tau) \rho$.
\end{lemma}
\begin{proof}
By \ref{lemma:J}, for every $g\in G$ 
$J=\det(\tau(g)) \det(\rho(g)) \rho(g) J$,
and hence $J\in E^*$.
\end{proof}
\begin{lemma}
\label{lemma:last}
Let $G$ be an extension of a planar symmetry group not of type R. 
Let $V\subset E$  denote the invariant plane,
Then $E^* \subset V$.
\end{lemma}
\begin{proof}
Let $V^*$ denote the 
orthogonal complement of the invariant plane.
Since $\det(\rho) = \det(\rho_2) \det(\rho_1)$ 
with $\det(\rho_2)\det(\tau) \neq 1$, the projection  
of $E^*$ on $V^*$ 
is fixed by the action of $G$ under the non-trivial homomorphism 
$\det(\tau)\det(\rho) \det(\rho_1) = \det(\tau) \det(\rho_2)$.
Hence $E^* \subset V$.
\end{proof}

By \ref{lemma:last}, one needs to consider the vectors
in the plane $V\subset E$  fixed by 
\[
\det(\rho_2(r))\epsilon_1 \rho_2(r)
\mbox{ \ and \ } 
-\det(\rho_2(h)) \epsilon_2 \rho_2(h)
\]
(the latter only if the action
type is not cyclic), 
where $\epsilon_1$ and $\epsilon_2$ are as above 
the elements $\rho_1(r)$ and $\rho_1(h)$.

For $C_6^+$, for example, $E^*$  is the subspace fixed by
$-\rho_2(r)$, which is a reflection along a line. Hence 
$C_6^+$, even if extension of a planar symmetry group not
of type R, might have minimizers with non-zero angular 
momentum. In fact, it is not difficult to see that 
up to a change in coordinates $C_6^+=C_3^-$.
On the other hand,
for $C_6^-$ it happens that $E^*$ is the subspace fixed by
$\rho_2(r)$, which is again a line. Again, as above,
$C_6^-$ can be written as $C_3^-$ with a suitable
choice of $\omega$ for $C_3^-$ (which is of type R).

Now we can consider the extensions of $D_{6}$ and $D_{12}$.
For $D_6$, we have  that
$\det(\rho_2(r))\epsilon_1 \rho_2(r)$
and $-\det(\rho_2(h)) \epsilon_2 \rho_2(h)$  
are respectively equal to
$\epsilon_1 $ 
and 
$\epsilon_2$
(seen as a $2\times 2$ matrices),
and hence equivariant minimizers of 
$D_{6}^{+,-}$, $D_{6}^{-,+}$
and $D_{6}^{-,-}$ have zero angular momentum.
On the other hand it is easy to see that after a change of coordinates 
$D_{6}^{+,+} = L_6^{+,-}$.
For $D_{12}$, 
$\det(\rho_2(r))\epsilon_1 \rho_2(r)$
and $-\det(\rho_2(h)) \epsilon_2 \rho_2(h)$  
are respectively equal to
$-\epsilon_1 \rho_2(r)$ (which is a reflection
along a line)
and 
$\epsilon_2$ (seen as a matrix).
Thus, if $\epsilon_2=-1$,
orbits have zero angular momentum. Otherwise, 
for $\epsilon_2=1$, it is not necessary. Furthermore, 
it is true that 
$D_{12}^{+,+} = L_6^{-,-}$
and $D_{12}^{-,+} = L_6^{+,-}$ for a suitable $\bfomega$.

The following definition is the natural extension
of the corresponding property for planar groups.
\begin{defi}
\label{defi:nottypeR3d}
A symmetry group $G$ is said \emph{of type R}
if there is a rotation axis for $G$ (and the restriction
of the action of $G$ on the invariant plane orthogonal
to the axis is a planar symmetry group of type R).
A symmetry group $G$ is said to be \emph{not of type R}
if there are not rotation axes for $G$ in $E$.
\end{defi}
The following lemma follows immediately from the 
previous arguments.
\begin{lemma}
\label{lemma:nottypeR}
If $G$ does not have rotation axes, i.e. 
it is not of type R,
than all $G$-equivariant trajectories have zero angular
momentum and hence they are planar.
\end{lemma}

\begin{lemma}
\label{lemma:secondkind}
Let $G$ be any space symmetry group  not of type R.
Then every 
$G$-equivariant non-collinear orbit is contained in a (unique)
$G$-invariant plane. 
\end{lemma}
\begin{proof}
By \ref{lemma:nottypeR}, the angular momentum is zero 
and hence the orbit is planar.
It is only left to show that the plane containing the 
orbit is $G$-invariant.
But since for every $g\in G$ and every $t\in \T$
\begin{eqnarray*}
\pm g \left[ (x_1(t) - x_2(t)) \times (x_1(t) - x_3(t)) \right] &= &
(gx_1(t) - gx_2(t))\times (gx_1(t) - gx_3(t)) \\
&=& 
(x_{g1}(gt) - x_{g2}(gt)) \times 
(x_{g1}(gt) - x_{g3}(gt),
\end{eqnarray*}
it follows that the plane containing the configuration
$x_1(t)$, $x_2(t)$ and $x_3(t)$ is $G$-invariant.
\end{proof}

\section{The vertical variation}
\label{sec:vertical}

Let $(z,w)$ be a system of coordinates for 
the Euclidean space $E\cong\R^3 \cong \CC \oplus \R$, 
with $z\in \CC$ and $w\in \R$. 
For a planar central configuration $\bar x$ consider
the  planar rotating periodic path $x(t)= e^{ikt} \bar x$,
with $k\in \ze$.  
In space 
the orbit can be written for $i=1,2,3$ 
as $(x_i(t),w_i(t))\in \R^2\times \R$ 
with $w_i=0$.
Now consider three periodic $H^1$-functions $\ff_i\from \T \to \R$.
There is a  corresponding path in $\Lambda$, which 
 will be denoted with
$(x(t),\varepsilon\varphi(t))$, obtained by adding the vertical
variation $\varepsilon\varphi$ to the rotation configuration $\bar x$.

\begin{lemma}
\label{lemma:secondd}
Let $\action(\varepsilon)$ denote the action of the path
$(x(t),\varepsilon\ff(t))$ in $[0,2\pi]$;
then
the second derivative of $\action(\varepsilon)$ evaluated
in $\varepsilon=0$ is
\[
\left.\dfrac{d^2\action}{d\varepsilon^2}\right|_{\varepsilon=0} =
\int_{0}^{2\pi} \left[
\sum_{i\in \nn} m_i \dot \ff_i^2
- \alpha
\sum_{i<j} \dfrac{m_im_j}{|x_i -x_j|^{\alpha+2} }
\left( \ff_i - \ff_j \right)^2
\right] \dt
\]
\end{lemma}
\begin{proof}
The second derivative of the kinetic part is
\[
\begin{split}
\dfrac{d^2}{d\varepsilon^2} \sum_{i\in\nn}
\dfrac{1}{2} m_i (|\dot x_i + \Omega x_i|^2 + \varepsilon^2 \dot \ff_i^2) =
\sum_{i\in \nn} m_i \dot \ff_i^2.
\end{split}
\]
Now, it is easy to see that
\[
\left.\dfrac{d^2}{d\varepsilon^2}\right|_{\varepsilon=0}
\left[
(a + \varepsilon^2b)^{c}
\right]
= 2 a^{c-1} b c;
\]
moreover, the terms in the potential part contain expressions of such type
\[
m_i m_i \left[ (x_i-x_j)^2  + \varepsilon^2(\ff_i-\ff_j)^2
\right]^{-\alpha/2},
\]
with $a= (x_i-x_j)^2$,
$b=(\ff_i - \ff_j)^2$ and
$c= -\alpha/2$.
Hence
\begin{eqnarray*}
\left.\dfrac{d^2}{d\varepsilon^2}\right|_{\varepsilon=0}
\sum_{i<j} m_im_j |(x_i,\ff_i) - (x_j,\ff_j)|^{-\alpha}
& = &  \sum_{i<j} m_im_j  2 \left[ (x_i-x_j)^2 \right]^{-\alpha/2-1}
(\ff_i - \ff_j)^2
(-\dfrac{\alpha}{2})
\\
&=&
-\alpha \sum_{i<j}  \dfrac{m_im_j}{ |x_i - x_j|^{\alpha+2}}
(\ff_i - \ff_j)^2. \\
\end{eqnarray*}
Thus, the claim.
\end{proof}
\begin{lemma}
\label{lemma:nicevar}
Consider the path 
$x(t)=e^{ikt}\bar x$  as above. For a unit 
vector $\ee\in \CC \subset  E$
define $\varphi=(\varphi_1,\varphi_2,\varphi_3)$ by the scalar
product
$\varphi_i(t) = x_i(t/k)\cdot \ee$ for $i=1,2,3$.
Then the second variation of \ref{lemma:secondd} is 
\[
\left.\dfrac{d^2\action}{d\varepsilon^2}\right|_{\varepsilon=0} =
\pi \left(  I(\bar x) - \alpha U(\bar x) \right),
\]
where $I(\bar x) = \sum_i m_i \bar x^2$ is the momentum
of inertia of $\bar x$ and 
$U(\bar x)$ the value of the potential function.
\end{lemma}
\begin{proof}
Define $\beta=1/k$; then 
\( \varphi_i(t) = x_i(\beta t) \cdot \ee 
=  (e^{it} \bar x_i) \cdot \ee
\) 
and therefore
\[
\int_0^{2\pi} m_i \dot \varphi_i^2 \dt =
\int_0^{2\pi} m_i \bar x^2 \sin^2(t + \delta_i) \dt  =
\pi m_i \bar x^2
\]
for some suitable $\delta_i$, which implies
that 
\[
\int_{0}^{2\pi} 
\sum_{i=1}^3 m_i \dot \ff_i^2 =
\pi \sum_{i=1}^3 m_i \bar x^2 = \pi I(\bar x).
\]
As for the second part of the expression in \ref{lemma:secondd},
since the norms $|x_i - x_j|$ are constant one obtains
\begin{eqnarray*}
\int_0^{2\pi} \dfrac{m_im_j}{ |x_i(t) - x_j(t)|^{\alpha+2}}
(\ff_i - \ff_j)^2  \dt &=&
\int_0^{2\pi} \dfrac{m_im_j}{ |x_i(t) - x_j(t)|^{\alpha}}
\left(\dfrac{\ff_i - \ff_j}{%
|x_i(t) - x_j(t)|
}\right)^2  \dt \\
&=&
\dfrac{m_im_j}{|\bar x_i - \bar x_j|^{\alpha}}
\int_0^{2\pi}
\left( \dfrac{x_i(\beta t) - x_j(\beta t)}{|x_i(\beta t) - x_j(\beta t)|}
\cdot \ee \right)^2 \dt \\
&=& 
\dfrac{m_im_j}{|\bar x_i - \bar x_j|^{\alpha}}
\int_0^{2\pi} \cos^2(t + \delta_{ij}) \dt \\
&=&
\pi \dfrac{m_im_j}{|\bar x_i - \bar x_j|^\alpha}.
\end{eqnarray*}
with a suitable choice of the shift constant $\delta_{ij}$.
Thus, summing up one obtains
\[
\int_{0}^{2\pi} 
\sum_{i<j} \dfrac{m_im_j}{|x_i -x_j|^{\alpha+2} }
\left( \ff_i - \ff_j \right)^2
\dt
= 
\pi
\sum_{i<j}
\dfrac{m_im_j}{|\bar x_i - \bar x_j|^\alpha} = \pi U(\bar x)
\]
The conclusion follows.
\end{proof}

Until now we did not assume anything else on $x(t)$ other 
than being rotating $k$ times during the interval $[0,2\pi]$.
Now we assume that it is a minimizer in a suitable 
\emph{linear} class of  paths (such as $\Lambda^G$ for a suitable
$G$ acting on the plane or the space).
Then the following equation holds
\begin{lemma}
\label{lemma:lastlem}
If $x(t)=e^{ikt}\bar x$ is  a minimizer of $\action_\omega$,
then  (Kepler's law)
\[
(k+\omega)^2 I(\bar x) = \alpha U(\bar x).
\]
\end{lemma}
\begin{proof}
It is easy to see that the action is 
\[
\dfrac{1}{2\pi}
\action_\omega = \dfrac{1}{2}(k+\omega)^2 I(\bar x) +  U(\bar x).
\]
By deriving the expression in $R=\sqrt{I}$ one obtains
that the minimum as $R>0$ varies is achieved 
when $(k+\omega)^2R^2 = \alpha U(\bar x)$. Otherwise,
one could also use homogeneity and directly
Newton's equations.
\end{proof}

\begin{propo}
\label{propo:vertical}
Assume that for a symmetry group $G$ every  rotating 
$G$-equivariant central configuration $x(t) = e^{ikt} \bar x$ is 
such that  $(k+\omega)^2>1$.
Then rotating central configurations 
cannot be minimizers of $\action^G$.
\end{propo}
\begin{proof}
By \ref{lemma:nicevar} the second variation 
is 
$\left.\dfrac{d^2\action}{d\varepsilon^2}\right|_{\varepsilon=0} =
\pi ( I(\bar x) - \alpha U(\bar x) )$.
But by  \ref{lemma:lastlem} 
\[
\alpha U (\bar x) =  (k+\omega)^2 I(\bar x),
\]
so that 
\[
\left.\dfrac{d^2\action}{d\varepsilon^2}\right|_{\varepsilon=0} =
\pi  I(\bar x) ( 1 - (k+\omega)^2 ) < 0,
\]
which shows that $x(t)$ cannot be a minimizer.
\end{proof}
Now we consider a different vertical variation,
which can be readily used for a vertical isosceles triangle.
Consider now the  variation $\varphi$ (of 
\ref{lemma:secondd}) defined as follows: $\varphi = v\sin t$,
where $v\in \R^n$ is a one-dimensional configuration
with $\sum_{i\in \nn} m_iv_i =0$. Without loss of generality
we assume that $\sum_{i\in \nn} m_i = 1$. 
\begin{lemma}\label{lem:othervar}
Let $G$ be a symmetry group and $x(t)=e^{ikt} \bar x$ a planar 
$G$-equivariant rotating central configuration.
If $(x,\varphi)$ is  $G$-equivariant and 
\begin{equation}\label{eq:othervar}
\sum_{i<j} m_i m_j(v_i - v_j)^2
\left(
1 - \alpha|\bar x_i - \bar x_j|^{\alpha+2}
\right)
<0,
\end{equation}
then 
$x(t)$ 
is  not a minimizer.
\end{lemma}
\begin{proof}
Since $\sum_{i\in \nn} m_i \varphi_i = 0$ and 
$\sum_{i\in\nn} m_i = 1$
by assumption,
one can write as in \ref{eq:differences}
 the kinetic energy in terms of differences,
and therefore equation \ref{lemma:secondd} can be read as
\[
\left.\dfrac{d^2\action}{d\varepsilon^2}\right|_{\varepsilon=0} =
\pi \left[
\sum_{i\in \nn} m_i m_j (v_i - v_j)^2
- \alpha
\sum_{i<j} \dfrac{m_im_j}{|x_i -x_j|^{\alpha+2} }
\left( v_i - v_j \right)^2
\right],
\]
since $\int_0^{2\pi} \sin^2 t = \int_0^{2\pi} \cos^2t = \pi$,
and this implies the claim.
\end{proof}
\begin{propo}
\label{propo:othervar2}
In the hypotheses of \ref{lem:othervar},
consider the following case: $n=3$, $m_1 = m_2$, 
$x_1 = - x_2$ (and hence $x_3=0$) and $v_1 = v_2$.
Then
$ \left.\dfrac{d^2\action}{d\varepsilon^2}\right|_{\varepsilon=0}  <0$
if and only if 
\begin{equation}
\label{eq:othervar2}
(k+\omega)^2 > \dfrac{m_1}{2^{\alpha+1}} + 1 - 2m_1.
\end{equation}
\end{propo}
\begin{proof}
If $x(t) = e^{ikt} \bar x$ 
is a minimum,
then 
$\bar x = (R,-R,0)$,
with 
\begin{equation}\label{eq:ralpha}
R^{\alpha+2} = \frac{\alpha}{2\beta},
\end{equation}
where 
\(
\beta = \dfrac{ (k+\omega)^2}{ m_1 2^{-\alpha} + 2m_3}\).
Since $v_1 = v_2$, the left term of \ref{eq:othervar} is
equal to 
\[
2  m_1 m_3(v_1 - v_3)^2
\left(
1 - \alpha|\bar x_1 - \bar x_3|^{\alpha+2}
\right),
\]
and therefore
$ \left.\dfrac{d^2\action}{d\varepsilon^2}\right|_{\varepsilon=0}  <0$
if and only if 
$R^{\alpha+2} < \alpha$.
But by \ref{eq:ralpha}
this is true if and only if 
\[
(k+\omega)^2 > \dfrac{m_1}{2^{\alpha+1}} + 1 - 2m_1,
\] 
as claimed.
\end{proof}

\begin{remark}\label{rem:anticip}
The right hand side of \ref{eq:othervar2} 
is a linear function of $m_1$, which is defined 
for $m_1 \in (0,1/2)$ and goes from a limit value of $1$ 
(for $m_1=0$) to a limit value of $2^{-(2+\alpha)}$ 
(for $m_1=1/2$).
Hence it is always possible to find  mass distributions
$(m_1,m_2,m_3)$ 
for which minimizers are not rotating Euler solutions, 
provided that for the minimal $k$ one has 
$(k+\omega)^2 > 2^{-(2+\alpha)}$.
For example, if $\alpha=1$ then one finds  that 
there is a non-trivial interval  of values of $\omega$
for which minimizers are non-trivial
for all  $m\in(\dfrac{3}{7},\dfrac{1}{2})$ (in the case  that
$k$ can be any integer) while the same happens
and for all $m\in(0,\dfrac{1}{2})$ 
if 
the symmetry group implies
a constraint on $k$ such that $k=0\mod 2$.
For equal masses  one has $m=\dfrac{1}{3} \not\in 
(\dfrac{3}{7},\dfrac{1}{2})$ and  so one has to assume
$k=0\mod 2$,  
inequality \ref{eq:othervar2} becomes
\[
(k+\omega)^2 > \dfrac{5}{12},
\]
and thus non-planar orbits exist for $\omega\in 
(\sqrt{\dfrac{5}{12}},
2 - \sqrt{\dfrac{5}{12}})$.
We will apply this simple argument below in \ref{lem:anticip} to prove
the existence of non-planar (quasi)-periodic orbits
for when two masses are approximately equal.
\end{remark}

\section{Space symmetries}
\label{sec:zzz}
In this section we will describe space symmetries and 
prove theorem \ref{theo:mt1}.

\begin{table}
\caption{The list of all dihedral symmetries}
\label{table:symmetries}
{\scriptsize
\begin{tabbing}
  \begin{tabular}{c}
  \(000\):
  \\%
  \(\left[
\begin{array}{cc}%
  +&+\\%
  +&+\\%
  +&+\\%
  \end{array}
\right]\)%
  \\\end{tabular}\=
  \begin{tabular}{c}
  \(001\):
  \\%
  \(\left[
\begin{array}{cc}%
  +&+\\%
  +&+\\%
  +&-\\%
  \end{array}
\right]\)%
  \\\end{tabular}\=
  \begin{tabular}{c}
  \(002\):
  \\%
  \(\left[
\begin{array}{cc}%
  +&+\\%
  +&+\\%
  -&+\\%
  \end{array}
\right]\)%
  \\\end{tabular}\=
  \begin{tabular}{c}
  \(003\):
  \\%
  \(\left[
\begin{array}{cc}%
  +&+\\%
  +&+\\%
  -&-\\%
  \end{array}
\right]\)%
  \\\end{tabular}\=
  \begin{tabular}{c}
  \(011\):
  \\%
  \(\left[
\begin{array}{cc}%
  +&+\\%
  +&-\\%
  +&-\\%
  \end{array}
\right]\)%
  \\\end{tabular}\\[18pt]
  \begin{tabular}{c}
  \(012\):
  \\%
  \(\left[
\begin{array}{cc}%
  +&+\\%
  +&-\\%
  -&+\\%
  \end{array}
\right]\)%
  \\\end{tabular}\=
  \begin{tabular}{c}
  \(013\):
  \\%
  \(\left[
\begin{array}{cc}%
  +&+\\%
  +&-\\%
  -&-\\%
  \end{array}
\right]\)%
  \\\end{tabular}\=
  \begin{tabular}{c}
  \(022\):
  \\%
  \(\left[
\begin{array}{cc}%
  +&+\\%
  -&+\\%
  -&+\\%
  \end{array}
\right]\)%
  \\\end{tabular}\=
  \begin{tabular}{c}
  \(023\):
  \\%
  \(\left[
\begin{array}{cc}%
  +&+\\%
  -&+\\%
  -&-\\%
  \end{array}
\right]\)%
  \\\end{tabular}\=
  \begin{tabular}{c}
  \(033\):
  \\%
  \(\left[
\begin{array}{cc}%
  +&+\\%
  -&-\\%
  -&-\\%
  \end{array}
\right]\)%
  \\\end{tabular}\\[18pt]
  \begin{tabular}{c}
  \(111\):
  \\%
  \(\left[
\begin{array}{cc}%
  +&-\\%
  +&-\\%
  +&-\\%
  \end{array}
\right]\)%
  \\\end{tabular}\=
  \begin{tabular}{c}
  \(112\):
  \\%
  \(\left[
\begin{array}{cc}%
  +&-\\%
  +&-\\%
  -&+\\%
  \end{array}
\right]\)%
  \\\end{tabular}\=
  \begin{tabular}{c}
  \(113\):
  \\%
  \(\left[
\begin{array}{cc}%
  +&-\\%
  +&-\\%
  -&-\\%
  \end{array}
\right]\)%
  \\\end{tabular}\=
  \begin{tabular}{c}
  \(122\):
  \\%
  \(\left[
\begin{array}{cc}%
  +&-\\%
  -&+\\%
  -&+\\%
  \end{array}
\right]\)%
  \\\end{tabular}\=
  \begin{tabular}{c}
  \(123\):
  \\%
  \(\left[
\begin{array}{cc}%
  +&-\\%
  -&+\\%
  -&-\\%
  \end{array}
\right]\)%
  \\\end{tabular}\\[18pt]
  \begin{tabular}{c}
  \(133\):
  \\%
  \(\left[
\begin{array}{cc}%
  +&-\\%
  -&-\\%
  -&-\\%
  \end{array}
\right]\)%
  \\\end{tabular}\=
  \begin{tabular}{c}
  \(222\):
  \\%
  \(\left[
\begin{array}{cc}%
  -&+\\%
  -&+\\%
  -&+\\%
  \end{array}
\right]\)%
  \\\end{tabular}\=
  \begin{tabular}{c}
  \(223\):
  \\%
  \(\left[
\begin{array}{cc}%
  -&+\\%
  -&+\\%
  -&-\\%
  \end{array}
\right]\)%
  \\\end{tabular}\=
  \begin{tabular}{c}
  \(233\):
  \\%
  \(\left[
\begin{array}{cc}%
  -&+\\%
  -&-\\%
  -&-\\%
  \end{array}
\right]\)%
  \\\end{tabular}\=
  \begin{tabular}{c}
  \(333\):
  \\%
  \(\left[
\begin{array}{cc}%
  -&-\\%
  -&-\\%
  -&-\\%
  \end{array}
\right]\)%
  \\\end{tabular}\end{tabbing}
}
\end{table}

\newcommand{\rS}{{r_\Sigma}}
\newcommand{\hS}{{h_\Sigma}}
\newcommand{\rV}{{r_V}}
\newcommand{\hV}{{h_V}}

Consider the case of groups with trivial core.
Let $r\subset G$ denote the $\T$-cyclic generator, and,
if it exists, let $h\in G$   denote one of the time reflections.
Consider $\rS=\sigma(r)$,
$\hS=\sigma(h)$, 
$\rV=\rho(r)$ 
and $\hV=\rho(h)$.
By \ref{lemma:sumofonedim}, the $G$-representation 
$\rho$ is the sum of one-dimensional components,
hence 
$\rV$ and $\hV$  can be written
as matrices with $\pm 1$ diagonal elements.
Thus a choice of the generators $r$ and $h$ yields
a $3\times 2$ matrix
\(\left[
\begin{array}{cc}%
  \rV^1&\hV^1\\%
  \rV^2&\hV^2\\%
  \rV^3&\hV^3\\%
  \end{array}
\right]\)
where the entries $\rV^i$  and $\hV^i$ are 
the diagonal entries of the matrices $\rV$ and $\hV$ respectively.
Conversely, if such a matrix is given, the elements
$r$ and $h$ can be obtained by the permutations $\rS$ 
and $\hS$ in $\Sigma_3$ (analogously for cyclic and brake action types).
The number of such matrices is the 
number of unordered $3$-tuples of elements 
chosen in the set $\{ [++],[+-],[-+],[--]\}$,
which are $\binom{4+3-1}{3} = 20$.
Under the identification $0=[++]$, $1=[+-]$,
$2=[-+]$ 
and $3=[--]$,
it is possible to represent such matrices with 
numbers of $3$ digits as in table
\ref{table:symmetries}.
 
If the group $G$  is of cyclic type, 
then there are  $12$ possible cases
for $\rV$ and $\rS$ ($4$ for $\rV$ times
$3$ for $\rS$): $\rS \in \{ (), (1,2), (1,2,3) \}$
and $\rV \in \{ [+++],[++-],[+--],[---]\}$.
Is not difficult to see that the resulting groups 
are three-dimensional
extensions listed in table \ref{table:cyclictype} 
(where, as we defined in \ref{defi:rotating}
 from now on the symbol `$\sim$' means 
that the symmetry group is equivalent to the group
in question after a change in rotating coordinates).
\begin{table}
\caption{Symmetry groups of cyclic type}
\label{table:cyclictype}\centering
\begin{tabular}{r|ccc}
 & $()$ & $(1,2)$ & $(1,2,3)$ \\
\hline
$[+++]$ & 
$C_1^+$ (fully uncoercive) & 
$C_2^+$ (fully uncoercive) & 
$C_3^+$ 
\\
$[++-]$ & 
$C_1^-$ &
$C_2^-$ (fully uncoercive) &
$C_3^-$
\\
$[+--]$ &
$L_2^-$ (fully uncoercive) &
$H_2^-$ (fully uncoercive) &
$C_6^- \sim C_3^+$
\\
$[---]$ & 
$\sim C_1^-$
& $\sim C_2^-$
& $\sim C_3^-$
\\
\end{tabular}
\end{table}
Given a symmetry group $G$, consider the elements
$\rV$, $\rS$, $\hV$ and $\hS$ defined above.
The matrix 
\(\left[
\begin{array}{cc}%
  \rV^1&\hV^1\\%
  \rV^2&\hV^2\\%
  \rV^3&\hV^3\\%
  \end{array}
\right]\)
associated to
$\rV$ and $\hV$ 
is one of the $20$ matrices of table \ref{table:symmetries}.
Furthermore, it is easy to prove that  up to 
permutations
the pair $[\rS,\hS]$ can be chosen
from the set
\[
\{
[ (1,2,3), (1,2) ],
[ (), (1,2) ],
[ (1,2), () ],
[ (1,2), (1,2) ],
[ (), () ]
\}.
\]
Consider first the case $[\rS,\hS] = [ (1,2,3), (1,2) ]$.
Since the matrices
\begin{equation}\label{eq:array}
\left[
\begin{array}{cc}%
  \rV^1&\hV^1\\%
  \rV^2&\hV^2\\%
  \rV^3&\hV^3\\%
  \end{array}
\right]
\mbox{ \ and \ }
\left[
\begin{array}{cc}%
  \rV^1&\hV^1\rV^1\\%
  \rV^2&\hV^2\rV^2\\%
  \rV^3&\hV^3\rV^3\\%
  \end{array}
\right]
\end{equation}
yield the same symmetry group up to change of coordinates
as  the matrix
it is possible to assume that $[\rV,\hV]$
belongs to one of the   $13$ items:
$000$, 
$001$, 
$002$ ($\cong 003$), 
$011$,
$012$ ($\cong 013$),
$022$ ($\cong 033$),
$023$,
$111$,
$112$ ($\cong 113$),
$122$ ($\cong 133$),
$123$,
$222$ ($\cong 333$),
$223$ ($\cong 233$).

Now, since we are excluding the case of 
groups which are bound to collisions, 
we have to rule out the cases in which 
$\hV$ or the product $\rV\hV$ is trivial, 
which means we have to exclude the  four 
cases 
$000$, 
$002\cong 003$,
$022 \cong 033$,
$222 \cong 333$; thus we are left with a list
of $9$ matrices.

We can now take the rotation axes into account.
Let $e_1$, $e_2$ and $e_3$ denote the canonical elements
of the base of the vector space $E$.
In table \ref{table:rotationaxes} it is found the list 
the elements of $\{e_1,e_2,e_3\}$ which are rotation
axes for the corresponding group.
\begin{table}
\caption{Rotation axes}
\label{table:rotationaxes}\centering
\begin{tabular}{|c|c|c|c|c|c|c|c|c|}
$001$ & 
$011$ & 
$111$ & 
$012$ & 
$112$ & 
$023$ & 
$122$ & 
$123$ & 
$223$ \\
\hline
 $e_1$, $e_2$ &
 $e_2$, $e_3$ &
 &
 $e_3$ &
 &
$e_1$ &
 &
 $e_1$ &
$e_2$, $e_3$\\
\end{tabular}
\end{table}
Using the rotating frame change of coordinates,
it is therefore possible to show 
that 
$023 \sim 001$,
$123 \sim 011$ and 
$223 \sim  012$, 
and hence 
that we are left with a choice among the $6$ items
$001=L_6^{+,-}$, 
$011=D_6^{+,+}\cong L_6^{+,-}$, 
$111 = D_6^{+,-}$, 
$012 = D_{12}^{+,+} \cong L_6^{-,+}$, 
$112 = D_6^{-,+}$ and 
$122 = D_{12}^{-,+}$.
Consider now the case $[\rS,\hS] = [ (), (1,2) ]$. 
We shall proceed in a similar way, analyzing case by case
until we are left with a small number of significative choices.
As above, $000$, $002$, $022$ and $222$ yield a group
which is bound to collisions, and  up to a change
of coordinates we can choose among the same following $9$
items:
$001$, 
$011$,
$012$ ($\cong 013$),
$023$,
$111$,
$112$ ($\cong 113$),
$122$ ($\cong 133$),
$123$,
$223$ ($\cong 233$).

The rotation axes are the same as 
those of table \ref{table:rotationaxes} above,
and again in a suitable rotating frame 
$023 \sim 001$,
$123 \sim 011$ and 
$223 \sim  012$, 
so that we can  
choose just among the $6$ items
$001$, $011$, $111$, $012$, $112$ and $122$.
First an easy computation shows that $111$, $112$ and $122$ 
yield symmetry groups 
without rotation axes and at the same time not coercive (thus
fully uncoercive by \ref{lem:pre}).
Furthermore,  the choice of $001$, $011$ or $012$ 
yields   groups $G$ conjugate to the 
three-dimensional 
extensions $H_2^{+,+}$,
$H_2^{+,-}$ and $H_2^{-,+}$
of the planar Isosceles  symmetry group $H_2$).
As a third possibility,  now consider the cases 
$[\rS,\hS] = [ (1,2), (1,2) ]$
and
$[\rS,\hS] = [ (1,2), () ]$.
We can consider just the case 
$[\rS,\hS] = [ (1,2), (1,2) ]$, up to a change of coordinates
(but we will not be able to use the argument  
of \ref{eq:array} to reduce the number of matrices).
As above, we start with the list \ref{table:symmetries} of all
possibilities. Since $\hS=(1,2)$, if $\hV$ is  trivial
then the resulting group is bound to collisions, 
and therefore we cancel the four  matrices $000$, $002$, $022$
and $222$. From the $16$ matrices left the following
$8$ do not have a rotation axis:
$003$,
$033$,
$111$,
$112$,
$113$,
$122$,
$133$ and
$333$.
Since the action on $\{1,2,3\}$ is not transitive,
by \ref{lem:pre}
all those with a row equal to $[++]$ are fully uncoercive
(that is, $003$ and $033$).
Furthermore, for the matrices
$111$, $112$ and $122$ 
the resulting symmetry group is bound to collisions 
(since $\rV\hV$ is the antipodal map,
while $\rS\hS$ is the trivial permutation).
About the three remaining items 
$113$, $133$ and $333$,
they are fully uncoercive simply because they contain
a row equal to $[--]$ (which yields a one-dimensional
non-coercive symmetry group).

So, we are left with $8$ choices, all with rotation axes:
$001$, 
$011$,
$012$,
$013$,
$023$,
$123$,
$223$,
$233$.
After a change in rotating coordinates
one can see that 
$023\sim 001$,
$123 \sim 011$,
$223\sim 012$
and $233\sim 013$; furthermore,
it is easy to see that $001$ and $013$ yield fully uncoercive
symmetry groups.
As a consequence, the remaining matrices
are $011$ and $012$.
In the notation of table \ref{table:thisclass},
they are respectively the symmetry groups
$H_4^{+,-}$ and
$H_4^{-,+}$.
At last, we can consider the case 
$[\rS,\hS] = [ (), () ]$ where 
the resulting group acts trivially on 
the index set.
The matrices 
$111$, 
$113$,
$133$
and $333$
yield bound to collisions groups.
By the same argument as \ref{eq:array},
we do not consider the duplicates 
$003$,
$013$,
$033$,
$112$,
$122$,
$222$,
$233$.
Of the remaining $9$ matrices,
three do not have rotation axes ($000$,
$002$ and $022$) and yield groups which are not coercive, 
while by a change in rotating
frames the other six can be reduced to the three
$001 (\sim 023)$, 
$011 (\sim 011)$ and 
$012 (\sim 012)$.
The group induced by $001$ is fully uncoercive.
So, we can  as before reduce the list of $20$ matrices to the two
cases $011$ and $012$ which yield
respectively the groups
$L_2^{+,-}$
and 
$L_2^{-,+}$ 
of table \ref{table:thisclass}.

\begin{defi}
\label{defi:vip}
We say that a symmetry group is of \emph{vertical isosceles} type
if its core  is generated by an element
$k$ such that $\rho(k)$ is the rotation of angle $\pi$
around an axis and $\sigma(k)$ is conjugate to 
the  permutation $(1,2)$.
\end{defi}

To conclude the proof it is left to prove the easy fact that 
if the core is non-trivial, then the group is homographic,
provided that it is not of  vertical isosceles type.

\begin{lemma}
\label{lemma:5.4}
Let $G$ be  a symmetry group which is not bound to collisions
and not fully uncoercive.
Then either $G$-equivariant trajectories are always homographic,
or the group is of vertical isosceles type
\ref{defi:vip}.
\end{lemma}
\begin{proof}
As in the proof of the (similar) proposition 5.4 for 
planar groups of \cite{zz}, 
$\ker\tau$ is isomorphic to a subgroup of $\Sigma_3$
and hence its three-dimensional representation  in $E$
is reducible: that is, there is a $\ker\tau$-invariant
line $\R$ in $E$. Now, since $\ker\tau$ is normal in $G$,
either $g\R \subset \R$ for every $g\in G$,
or $E$ is the sum of three copies of $\R$, 
which are permuted by the elements in $G$ (hence $\ker\tau$
has order $2$ and acts as the antipodal map $-1$ on $E$).
After an analysis of a few cases, it is clear that 
the possible actions of  the core are the following:
$\ker\tau = \left< (a_3 ,(1,2)) \right> $,
$\ker\tau = \left< (r_2 ,(1,2)) \right> $,
$\ker\tau = \left< (r_3 ,(1,2,3)) \right> $,
$\ker\tau = \left< (r_3,(1,2,3),(h_2,(1,2)) \right>$
(which, incidentally, are extensions of 
planar analogues), where 
$a_3$ is the antipodal map $a_3=-1$,
$r_2$ the rotation of $\pi$ around a fixed axis,
$r_3$ the rotation of $2\pi/3$ around a fixed axis $\R$,
and $h_2$ the reflection with respect to a plane
containing (or with respect to a line orthogonal to $\R$).
In the first case $\ker\tau$-invariant configurations
are antipodal binaries with a third mass at 
the origin (this is equivalent to a spatial Kepler problem),
which is a homographic group;
the second is the well-known case  of isosceles triangle (which is
not homographic);
the third and the fourth yield homographic symmetry groups.
Thus the proof.
\end{proof}
\section{Proof of theorem \ref{theo:mt2}}
\label{sec:main}
In this section we will complete the proof of \ref{theo:mt2},
proving one-by-one  all its parts.
\begin{lemma}
\label{lemma:collisions}
Let $G$ be a symmetry group of the three-body problem
which is not bound to collisions. Then 
all $G$-equivariant minimizers  are collisionless.
\end{lemma}
\begin{proof}
Consider first the case of a group with trivial core; let $H$ be 
one of its maximal $\T$-isotropy groups: $H$ is generated 
by the non-trivial element $h\in H$ with
$\sigma(h) \in \{ (), (1,2) \}$.
The orthogonal motion
$\rho(h)$ is of order at most two, and hence 
there are at least three invariant orthogonal lines (equivalently,
the representation of $H$ is the sum of one-dimensional
$H$-representations). An immediate consequence is 
that if it is not bound to collisions
(that is, if $(\rho(h),\sigma(h)) \neq (1,(1,2))$, 
the subgroup $H$ has the rotating circle property  (10.1)
of \cite{FT2003}. Thus, by (10.10) of the same paper, 
$G$-equivariant minimizers
are collisionless.  
Now assume that the kernel of $\tau$ is not trivial.
As in the proof of \ref{lemma:5.4} it is possible
to assume that $\ker\tau$ acts as one of the four 
cases listed, where the first
case yields a one-center (Kepler) problem in space,
the third  and fourth cases yield a one-center planar problem,
and 
the second case is the isosceles triangle,  see \ref{defi:vip}.
One can readily see that theorem (10.1) of \cite{FT2003}
can be applied  in all the cases, and hence the thesis.
\end{proof}
\begin{lemma}
\label{lemma:mt2pre1}
For every $\omega$, minimizers for 
$C_1^-$,
$H_2^{+,-}$,
$C_3^+$,
$L_6^{+,+}$,
$L_6^{+,-}$
are Lagrange homographic solutions.
\end{lemma}
\begin{proof}
It is easy to see that the hypothesis of \ref{propo:gordon} 
are fulfilled, hence the thesis.
\end{proof}

\begin{lemma}\label{lemma:mt2pre2}
Minimizers for $D_6^{+,-}$,
$D_{6}^{-,+}$
and $D_{12}^{-,+}$
are zero-angular momentum planar solutions: the
Chenciner-Montgomery figure-eight solution.
\end{lemma}
\begin{proof}
Since these groups do not have rotation axes,
by lemma ref{lemma:nottypeR}, their minimizers 
are planar orbits with zero angular momentum,
contained in a $G$-invariant plane by \ref{lemma:secondkind}.
All planes are $G$-invariant for $D_6^{+,-}$,
and the restricted group is $D_6$. Hence  
by (4.15) of \cite{zz}
the minimizer for $D_6^{+,-}$ 
is the $D_{12}$-symmetric Chenciner-Montgomery figure eight \cite{monchen}.
Next, if $G=D_{12}^{-,+}$,
a $G$-orbit (which is collisionless by \ref{lemma:collisions},
has to be contained in the
$G$-invariant plane 
where $G$ acts as $D_{12}$, since otherwise 
it would result to be bound to collisions. Hence the minimum 
for $D_{12}^{-,+}$ is again Chenciner-Montgomery figure eight.
At last, consider $D_{6}^{-,+}$, which is a group of order $12$.
In any one of the infinitely many invariant planes
with $D_{12}$-action as restriction there is a CM-eight minimum,
while in the other invariant plane there is a
redundant $D_6$-action (hence if the minimum were to  be
contained in this plane, it would have implied that 
the action of the $D_6$-eight is less than the action
of the $D_{12}$-eight, which is not true by \cite{zz}).
Thus as above the minimum is the $D_{12}$-symmetric
CM-eight.
\end{proof}

\begin{lemma}
\label{lemma:mt2pre3}
For $\omega\in(-1,1) + 6\ze$,  minimizers for 
$L_6^{-,+}$ and its subgroup 
$C_3^- 
\subset L_6^{-,+}$ are the non-planar (if $\omega\neq 0$)  families of 
quasi-periodic solutions
called respectively $P_{12}$ and $P'_{12}$, which might be
likely to coincide.
\end{lemma}
\begin{proof}
By \ref{lemma:collisions}, 
minima for $C_3^-$ and $L_6^{-,+}$ 
exist and are collisionless.
If they are planar, then in both cases
they have to be a Lagrange rotating solution
(rotating at angular speed $k$, with
$k=\pm 2 \mod 6$, minimizing the number $(k+\omega)^2$).
The proof can be concluded by applying \ref{propo:vertical}.
Action levels for the resulting minima
are depicted in figure \ref{fig:p12}: in the intervals $(1/6,1/2)$
and $(1/2,5/6)$  the minimum is a rotating Lagrange 
solution, while otherwise it is a non-planar
orbit
(in the graph the period is rescaled to $12\pi$,
so that values of $\omega$ need to be multiplied by
a factor $6$). It is remarkable how the estimate
of \ref{propo:vertical}  seems to be sharp, since
when its hypothesis is not fulfilled  
one can find numerically that minimizers are in fact planar rotating
Lagrange triangles. 

\begin{figure}\centering
\caption{Action-levels for the $L_{6}^{-,+}$-symmetric
$P_{12}$ minimizers}
\label{fig:p12}
\psfrag{action}{$\action_\omega$}
\psfrag{omega}{$\omega$}
\includegraphics[width=0.7\textwidth]{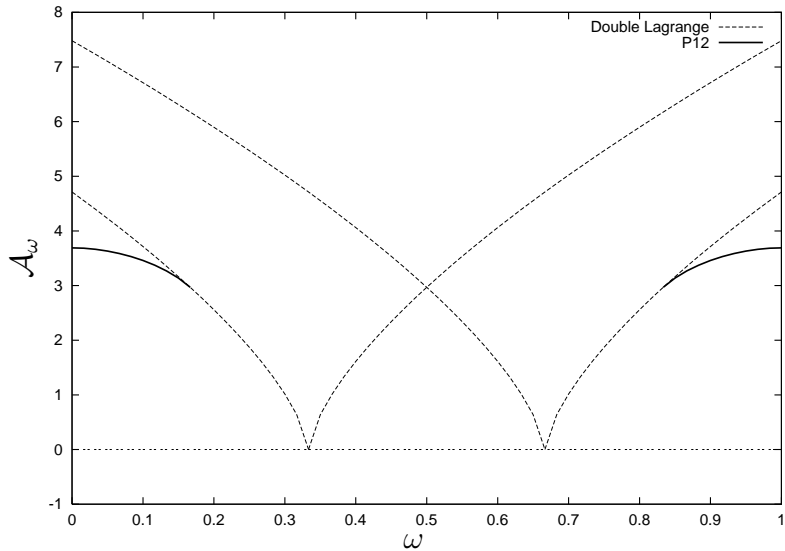}
\end{figure}
\end{proof}

\begin{remark}
\label{rem:p12}
In \cite{marchalp12}\footnote{%
The proof of the existence, claimed in \cite{marchalp12}
and later in \cite{marchal}, was completed in
\cite{chenciner_taiyuan,chenICM}   
by action level estimates on colliding paths,
since the main result of \cite{marchal} cannot be
applied to the $P_{12}$ family, which has 
a symmetry group of dihedral type. The action levels graph
of \cite{chenciner_taiyuan} is a qualitative picture
of figure \ref{fig:p12}}
Marchal actually introduced, under a different notation, 
the family corresponding to minimizers 
(as $\omega$ varies) under  $L_{6}^{-,+}$-symmetries,
naming it $P_{12}$  -- see also \cite{marchal_book}.
The fact that there is also a $C_3^-$-equivariant 
family (as well as the fact that there is a
group action of cyclic type yielding a figure-eight orbit)
apparently was not known. 
As for the planar eights, it seems that (numerically) these
two families coincide, and that 
for $\omega=0$ one finds the CM-eight and the
cyclic eight (which should as well likely coincide).
Some of these 
were questions raised in the last section 
of \cite{marchalp12} (see also section 4.(iv) of \cite{chenICM}),
questions which  probably
still need to find an answer.  
For example, as said at the end of section 4 of \cite{chenICM},
one should prove that  for $\omega=0$ the minimum is planar
and that minimizers are a continuous family 
(from figures \ref{fig:p12}  and \ref{fig:p12paths} it seems that 
action levels depend continuously on
$\omega$, as well as the corresponding trajectories). 
In figure \ref{fig:p12paths} 
the minimizers corresponding to the 
values $\omega=j/5$, for $j=0\dots 5$ are shown,
together with their projections. The curves with $j=0$ and 
$j=5$ have a bigger width. On its side,
there is one of the orbits (actually, corresponding
to the non-integer value $j=2.5$) in the inertial frame.
\begin{figure}\centering
\subfigure[Some representatives of the 
$P_{12}$ family in the rotating frame]{%
\includegraphics[width=0.45\textwidth]{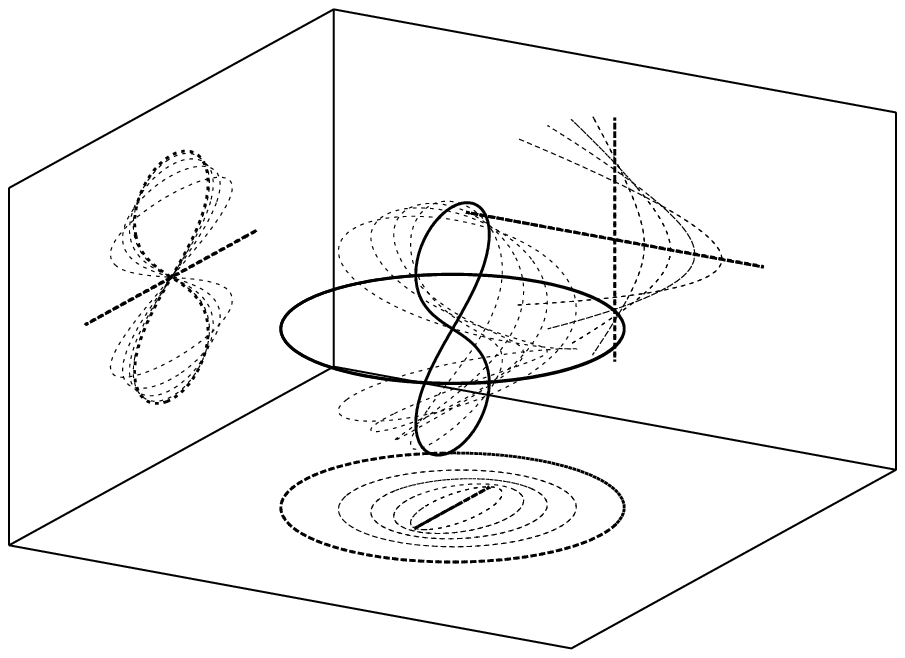}
\label{fig:p12paths}}
\hspace{0.04\textwidth}
\subfigure[One of representatives of the $P_{12}$ family on inertial frame]{%
\includegraphics[width=0.45\textwidth]{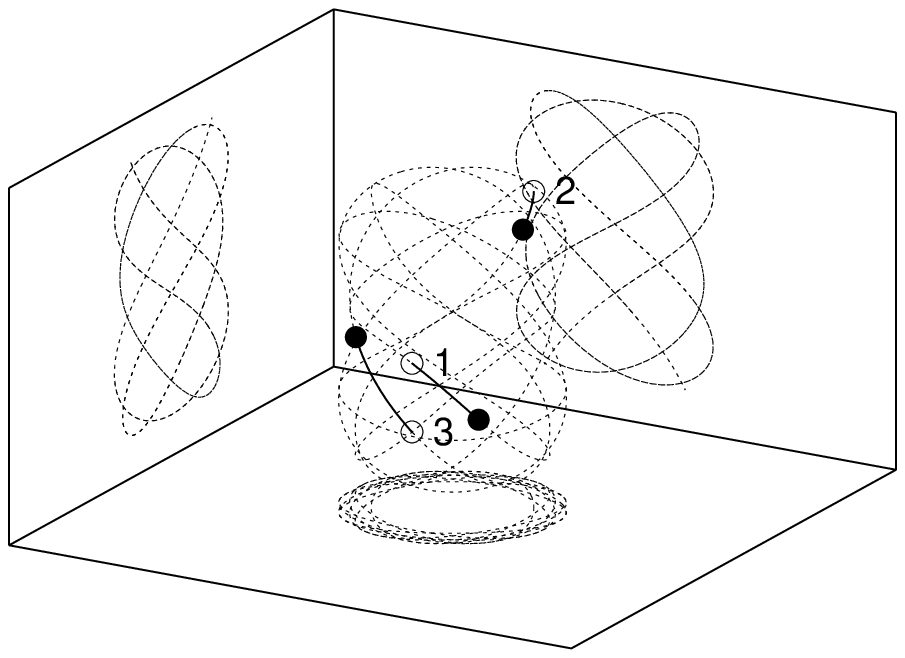}
\label{fig:p12if}}
\caption{The $P_{12}$ family}
\end{figure}
\end{remark}

\begin{remark}
In \cite{cfm}, Chenciner, F\'ejoz and Montgomery
found, under an assumption (at the moment numerically evaluated)
of non-degeneracy for the CM-eight,
three families of periodic orbits in rotating frames,
originated by three different symmetry-breaks of the planar
eight. 
The family there termed $\Gamma_1$ is $P_{12}$. 
\end{remark}

\begin{lemma}\label{lem:anticip}
For every $\omega\not\in \ze$ there is a mass distribution
such that minimizers of $L_2^{+,-}$,
$L_2^{-,+}$, 
$H_4^{+,-}$ and 
$H_4^{-,+}$ (where for the first two masses do not need to be equal,
while in the second three two masses need to be equal) 
are not planar.
Conversely,
for every mass distribution with two equal masses there is an $\omega$
such that 
minimizers symmetric under the groups
$L_2^{-,+}$ are
not planar.
\end{lemma}
\begin{proof}
If follows directly from remark \ref{rem:anticip}.
\end{proof}
\section{Remarks}
\label{sec:remarks}

We did not prove general and complete results about planarity or non-planarity
of minimizers for the following 
symmetry groups:
$L_2^{+,-}$,
$L_2^{-,+}$,
$H_2^{-,+}$,
$H_4^{+,-}$
and $H_4^{-,+}$. In fact, we have proved that  for some
choices of masses and angular speed minima are non-planar,
but we could not prove that minima are planar for other choices
(as apparently they are: all the Hill-type orbits
and Euler solutions exposed in for example \cite{zz}).
We describe now, among other remarks,  
some properties of their minimizers
with a little bit more details.

\begin{rem}
If the masses are equal,
for all $\omega$ 
minimizers under the symmetry group $L_2^{+,-}$ 
are planar (and they are the Euler and Hill retrograde orbits
described in \cite{zz} --
note that 
it is easy to prove by \ref{propo:othervar2} 
that for all these symmetry groups 
there are choices of masses and angular speed 
$\omega$ for which 
minimizers are not planar).
Also, it turns out that there are other local minimizers (planar 
and non-planar).
\end{rem}
\begin{rem}
The symmetry group $L_2^{-,+}$ imposes that possible
rotating configurations (which are Euler collinear) 
have to rotate an even number of loops, i.e.
the cyclic part of the symmetry group
imposes that a rotating central configuration has 
to have $k = 0 \mod 2$, 
and hence by \ref{rem:anticip} there is a continuum of choices
for $\omega$, for every choice of non-zero masses,
such that the minimizer is non-planar.
An example is shown in figure \ref{fig:nonplanar}.
\begin{figure}\centering
\subfigure[Rotating frame]{%
\includegraphics[width=0.45\textwidth]{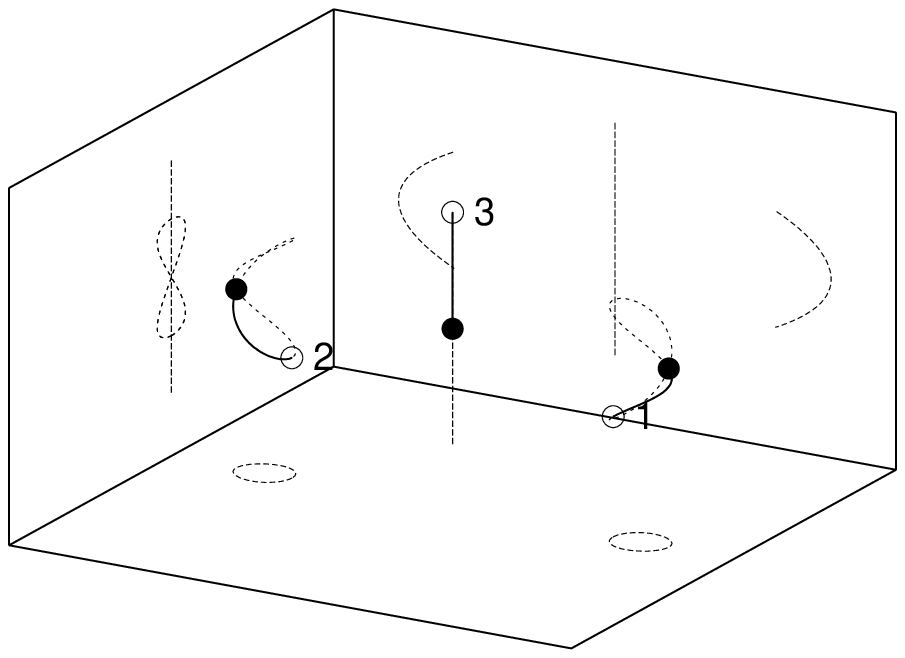}
\label{fig:nonplanar_rot}
}
\hspace{0.05\textwidth}
\subfigure[Inertial frame]{%
\includegraphics[width=0.45\textwidth]{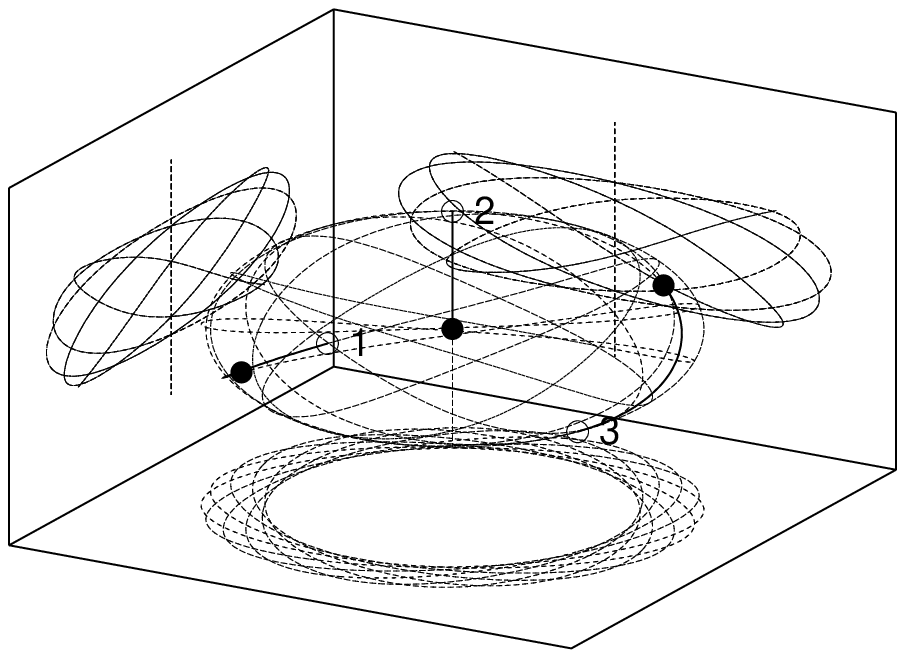}
\label{fig:nonplanar_if}
}
\caption{A non-planar $L_2^{-,+}$-symmetric orbit with equal masses, $\omega=0.9$}
\label{fig:nonplanar}
\end{figure}
\end{rem}
\begin{rem}
If the symmetry group is $H_2^{-,+}$, then  again it implies that 
a rotating central configuration has 
to rotate a number $k = 0 \mod 2$ of times. The possible
rotating configuration
is a Lagrange   triangle and it is not possible to
apply \ref{propo:vertical} to show that it is not a minimum
since  it is always possible
to find $k$ with 
$|k+\omega|\leq 1$. On the other hand the constraint
$k=0\mod 2$ prevents one to apply \ref{propo:gordon}
to show that the minimum is in fact
a Lagrange solution. Numerical experiments
show that this is the case.
\end{rem}
\begin{rem}
Now consider the symmetry groups 
$H_4^{+,-}$
and $H_4^{-,+}$.
As in the planar case, 
one finds non-homographic minimizers for 
some interval of 
values of $\omega$ (and approximately equal masses  -- see figure
4 of \cite{zz}). Homographic solutions have to 
be Euler-Moulton rotating collinear configurations,
and hence we can literally repeat the arguments
applied above to  $L_2^{+,-}$
and $L_2^{-,+}$,
and obtain the fact that for two equal masses 
it is always possible to find intervals 
of angular speed $\omega$ 
for which minimizers are not planar. 
This time rotating central configurations 
do not necessarily have $k=0 \mod 2$,
since $H_4$ has already a cyclic part of order $2$,
and hence one cannot prove with the vertical variation 
above that there are 
non-planar orbit for equal masses 
(for equal masses, after numerical experiments
it seems that local minima under  $H_4^{+,-}$ 
are planar,
while local minima under the action of $H_4^{-,+}$
can be non-planar for some $\omega$).
\end{rem}

\begin{rem}
\label{rem:extendit}
In the proofs of \ref{propo:vertical} 
and \ref{lem:othervar}
there is clearly 
no need to assume the bodies to be three.
In fact, the same vertical variation yield 
interesting non-planar orbits for every number
of bodies. 
For example,  it is easy to 
show by a straightforward extension of \ref{lemma:nicevar}
that families corresponding to the $P_{12}$ and $P'_{12}$ ones 
exist for any
odd number $n$ of bodies. In fact, consider the cyclic group $C$
of order $2n$, acting by a cyclic permutation of the $n$ bodies
on the index set and by a reflection along a plane $p$ in the space $E$.
Then, $C$-symmetric loops are choreographies in $E$
consisting in $n$ bodies, and if we choose as rotation axis $\bfomega$
the line orthogonal to the plane $p$,
we obtain family of coercive $n$-body problems such 
that for  $\omega \in (-1,1) \mod 2n$
minima are not equilibrium solutions. Since they are
collisionless due to \cite{FT2003},
they are periodic orbits (non-planar for $\omega\neq 0$).
It is likely that they behave like the $P_{12}$, namely
that they connect a eight solution with a  (twice rotating in the 
rotating frame) homographic
solution, but this is probably hard to prove 
(already for $n=3$ nobody has published a proof 
that for $\omega=0$ the $P_{12}$ is a planar eight
and that the family is a continuous one).
\end{rem}

\begin{rem}
About \ref{propo:vertical}, 
a similar proposition was  used by Chenciner 
in \cite{chenkyoto} to show that 
minimizers for the $n\geq 4$ anti-symmetric loops
are non-planar, following Moeckel's theorem
on central configurations \cite{Mo90}.
While the computation is very similar, here 
we use a different type of vertical variation,
which yield solutions in particular also when the rotating 
central configuration minimizes the reduced potential $\tilde U$,
due to the greater number of loops that symmetry 
constraints impose on rotating central configurations.
Furthermore, in section 3 of \cite{chenkyoto} 
there is an interesting 
short  remark proving the existence of a  non-planar periodic solution
for 3 bodies in the vertical isosceles problem 
under the antisymmetry constraint.
Since this constraint coincides with  
the group $C_1^-$ with angular speed $\omega=1$ (group 
which implies $k=0\mod 2$ to any $C_1^-$-equivariant equilibrium solution),
one can use \ref{propo:othervar2} 
and obtain that for all 
$\omega\in (\sqrt{\dfrac{5}{12}},2 - \sqrt{\dfrac{5}{12}})$ 
(as in remark \ref{rem:anticip} above)
and equal masses
a $C_1^-$-symmetric vertical isosceles minimizer 
is not planar. In figure \ref{fig:isosce}
it is shown the solution in the inertial frame --
probably this is the simplest non-planar periodic solution
of the three-body problem.
\begin{figure}
\centering
\includegraphics[width=0.4\textwidth]{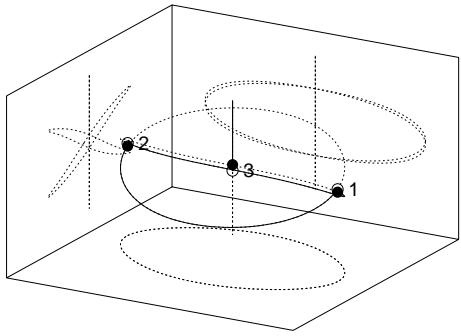}
\caption{The simplest non-planar solution with three equal masses}
\label{fig:isosce}
\end{figure}
\end{rem}

\begin{rem}
\label{rem:omegas}
In this paper we have sometimes rescaled the period
to a number different than $2\pi$. 
It is worth mentioning that, because of the homogeneity of the potential,
a  minimizer in a frame rotating with  angular velocity $\omega$ and
period $k 2\pi$ 
corresponds to a minimizer with period $2\pi$ in a frame rotating with 
angular velocity $k\omega$. The reason of rescaling the period
is that in our numerical experiments we have 
decided to rescale the period in order to have 
a fundamental domain of length $\pi$. 
The data for the minimizers used for figures were
obtained by a custom optimization program running on a Linux
cluster. The symmetries are computed by a package written in  GAP 
and python.
\end{rem}

\begin{rem}
In theorem \ref{theo:mt2}
nothing is stated about $H_2^{-,+}$.
Numerically one finds that $H_2^{-,+}$-symmetric minima
are rotating Lagrange triangles,
but \ref{propo:gordon} blow  cannot be
applied to this case. 
\end{rem}

\def\cfudot#1{\ifmmode\setbox7\hbox{$\accent"5E#1$}\else
  \setbox7\hbox{\accent"5E#1}\penalty 10000\relax\fi\raise 1\ht7
  \hbox{\raise.1ex\hbox to 1\wd7{\hss.\hss}}\penalty 10000 \hskip-1\wd7\penalty
  10000\box7} \def\cprime{$'$} \def\cprime{$'$} \def\cprime{$'$}
  \def\cprime{$'$} \def\cprime{$'$} \def\cprime{$'$} \def\cprime{$'$}
  \def\cprime{$'$} \def\cprime{$'$} \def\cprime{$'$} \def\cprime{$'$}
  \def\cprime{$'$} \def\cprime{$'$} \def\cprime{$'$} \def\cprime{$'$}
  \def\cprime{$'$} \def\polhk#1{\setbox0=\hbox{#1}{\ooalign{\hidewidth
  \lower1.5ex\hbox{`}\hidewidth\crcr\unhbox0}}}
  \def\polhk#1{\setbox0=\hbox{#1}{\ooalign{\hidewidth
  \lower1.5ex\hbox{`}\hidewidth\crcr\unhbox0}}} \def\cprime{$'$}
  \def\cprime{$'$} \def\cprime{$'$}

\end{document}